\documentclass[12pt]{amsart}
\usepackage[
  dvipdfmx,
  colorlinks=true,
  linkcolor=black,
  citecolor=black,
  urlcolor=blue
]{hyperref}
\usepackage{color}
\usepackage{amsmath,amssymb,amsthm}
\usepackage{amscd}
\usepackage{comment}
\usepackage{overpic}
\usepackage{tabularx}
\usepackage{array}
\usepackage[dvipdfmx,left=25truemm,right=25truemm,top=30truemm,bottom=30truemm]{geometry}

\usepackage{bm}

\newtheorem{theorem}{Theorem}[section]

\newtheorem{proposition}[theorem]{Proposition}

\theoremstyle{definition} 
\newtheorem{remark}[theorem]{Remark}

\newcommand{\R}{\mathbb{R}}

\renewcommand{\epsilon}{\varepsilon}

\renewcommand{\geq}{\geqslant}
\renewcommand{\leq}{\leqslant}

\makeatletter
\@addtoreset{equation}{section}

\makeatother

\begin{document}

\title[On the contact of regular surfaces in $\R^4 with toroidal cylinders$]
{On the contact of regular surfaces in $\R^4$ with toroidal cylinders and their apparent contours}

\author[M.~Hasegawa]
{Masaru Hasegawa}

\address[Masaru Hasegawa]{%
Department of Information Science, Center for Liberal Arts and Sciences, Iwate Medical University,
1-1-1 Idaidori, Yahaba-cho, Shiwa-gun, Iwate 028-3694, Japan.}
\email{mhase@iwate-med.ac.jp}


%


\allowdisplaybreaks[4]


\maketitle


\begin{abstract}
We study the contact of regular surfaces in $\mathbb{R}^4$ with toroidal  cylinders $S^1(r) \times \mathbb{R}^2$. 
We relate the contact type to the singularities of the corresponding parallel projections along the direction of the axes of the toroidal cylinders, the geometry of the apparent contours under these projections, and the geometry of the original surfaces.
\end{abstract}


\section{Introduction}
\label{sec:introduction}

In this paper, we characterize the contact type between regular surfaces in $\R^4$ and toroidal cylinders $S^1(r) \times \R^2$, which are one of homogeneous surfaces in $\R^4$, in terms of information obtained from the corresponding parallel projections of the surfaces onto planes, the geometry of their apparent contours and the geometry of the original surfaces.

The study of contact between submanifolds and model submanifolds (for example, lines, circles, planes, spheres, etc.) from the viewpoint of singularity theory goes back to Thom's work \cite{Thom1972}. 
Following Thom's suggestion, Porteous \cite{Porteous1971} studied contact with hyperspheres using distance-squared functions, while Montaldi \cite{Montaldi1986-1} subsequently established a general framework for contact between submanifolds in terms of Mather's framework \cite{Mather1968} of $\mathcal{K}$-equivalence.


The geometry of regular surfaces in $\R^4$ has been studied from the viewpoint of the above framework of contact and singularity theory.
Mochida, Romero-Fuster and Ruas \cite{MRFR1995} studied contact with hyperplanes, whereas Romero-Fuster and Sanabria-Codesal \cite{RFSC2013} studied contact with hyperspheres. 
More recently, the author \cite{Hasegawa2026-2} studied contact with spherical cylinders, which are another type of homogeneous surfaces in $\R^4$.  

Toroidal cylinder is a member of homogeneous surfaces in $\R^4$. 
A homogeneous surface in Euclidean $n$-space $\R^n$ is an orbit of a certain subgroup of the Euclidean motion group $O(n) \rtimes \R^n$.
This notion provided by Takahashi \cite{Takahashi1970} (see also \cite{Segre1938}). 
To explain it in terms of the extrinsic geometry, the homogeneous surfaces have constant principal curvatures at every point. 
In $\R^n$, the homogeneous surfaces are hyperplanes $\R^{n-1}$, hyperspheres $S^{n-1}(r)$ and hypercylinders $S^k(r) \times \R^{n-k-1}$.  
In particular, in $\R^3$, there are only three types homogeneous surfaces: planes $\R^2$, spheres $S^2(r)$ or cylinders $S^1(r) \times \R$.  
Moreover, in $\R^4$, there are only four types: hyperplanes $\R^3$, hyperspheres $S^3(r)$, spherical cylinders $S^2(r) \times \R$ and toroidal cylinders $S^1(r) \times \R^2$. 

The contact of surfaces, including singular surfaces, in $\R^3$ with planes or spheres has been widely studied (see, for example,  \cite{AOS2026, BGT1995, BGT1999, Francisco2023, FH2012-1, FH2012-2, FH2013, FH2025-1, FH2025-2, OSS2018, OST2018, Porteous1983, Teramoto2026}). 
On the other hand, the contact between regular surfaces in $\R^3$ and cylinders had not been investigated in detail until the work of Fukui, the author and Nakagawa \cite{FHN2017}.   
In \cite{FHN2017}, they introduced a new geometric notion called a cylindrical direction of the surfaces, which arises from contact with cylinders. 
A cylindrical direction $\bm{v}\in S^2$ at a point $\bm{p}$ of a regular surface is the direction of the axis of a cylinder having degenerate contact with the surface, where the axis passes through the finite focal point of its regular apparent contour along $\bm{v}$, and the contour has a vertex at the point corresponding to $\bm{p}$.
Tari, Salarinoghabi and the author \cite{TSH2025} investigated  deformations of apparent contours of regular surfaces in $\R^3$, using the fact that cylindrical directions capture vertices of their apparent contours. 
The author and Saji \cite{HS2025} investigated whether the notion of cylindrical directions could be extended to cross caps and cuspidal edges through the existence of cylinders having degenerate contact with these surfaces. 
They showed that an analogue of cylindrical directions exists for cross caps, whereas the notion cannot be extended to cuspidal edges. 
The author \cite{Hasegawa2026-2} introduced a higher-dimensional analogue, called a spherical-cylindrical direction, for regular surfaces in $\R^4$. 
This direction captures the ridge condition of a projected regular surface under an orthogonal projection along the axis of a spherical cylinder having degenerate contact with the original surface.

The contact of regular surfaces in $\R^4$ with toroidal cylinders is related to parallel projections along the axes of the cylinders. 
Nogueira \cite{Nogueira1998} investigated the relationships among the bifurcation sets of parallel projections of surfaces in $\R^4$ onto planes, height functions, and contact with 2-planes, which had not been studied in detail before his work. 
Deolindo-Silva and Salarinoghabi \cite{DSS2026} showed that contact with 2-planes can be characterized by singularities of parallel projections along these planes onto transverse planes. 
They also investigated inflections and vertices of the apparent contours of the surfaces under these projections.

Our main results are described in Theorems~\ref{thm:A2} and \ref{thm:Ak}. 
Theorem~\ref{thm:A2} characterizes the contact type between a regular surface in $\R^4$ and a toroidal cylinder according to the rank of the Jacobian matrix of the corresponding parallel projection. 
Theorem~\ref{thm:Ak} shows that, under the fold condition in Theorem~\ref{thm:A2} (1), which ensures that the apparent contour is regular and non-inflectional, the contact type is characterized by the order of the vertex of the apparent contour.
This provides a higher-dimensional analogue of the cylindrical direction of a regular surface in $\R^3$, given by the direction of the axis of the toroidal cylinder.

The paper is organized as follows.
In Section~\ref{sec:Preliminaries}, we recall the differential geometry of regular surfaces in $\R^4$, the notion of contact from the viewpoint of singularity theory, and basic results concerning cylindrical directions of regular surfaces in $\R^3$. 
In Section~\ref{sec:contact_with_toroidal_cylinders}, we characterize the contact type between regular surfaces in $\R^4$ and toroidal cylinders in terms of the singularities of the corresponding parallel projections, the geometry of their apparent contours, and the geometry of the original surfaces.


\section{Preliminaries}
\label{sec:Preliminaries}

\subsection{Differential geometry of surfaces in $\R^4$}
\label{sec:DGS}

The differential geometry of surfaces in $\R^4$ has been well studied from various perspectives (see, for example, \cite{BN1998, BT2002, GMRFR2000, Little1969, MRFR1995, NBT2007}). 
In particular, \cite[Chapter 7]{IRFRT2015} provides a systematic overview of this topic from the viewpoint of singularity theory. 

Let $U$ be an open subset of $\R^2$, and let a regular surface $M \subset \R^4$ be locally parameterized by
$$
X\colon U \to \R^4
$$
such that $X(\bm{q})=\bm{p}\in M$. 
Then $\{X_u(\bm{q}),X_v(\bm{q})\}$ is a tangent frame of the tangent space $T_{\bm{p}}M$. 
The \textit{first fundamental form at $\bm{p}$}
$$
I_{\bm{p}}\colon T_{\bm{p}}M \times T_{\bm{p}}M \to \R
$$
is defined by
$$
I_{\bm{p}} = \langle dX, dX\rangle = E du^2+ 2F du dv+G dv^2,
$$
where
$$
E=\langle X_u, X_u\rangle, \quad F=\langle X_u, X_v\rangle, \quad G=\langle X_v, X_v \rangle.
$$
Let $\{\bm{\nu}_1,\bm{\nu}_2\}$ be an orthonormal frame of the normal space $N_{\bm{p}}M$. 
The \textit{second fundamental form at $\bm{p}$}
$$
II_{\bm{p}}\colon T_{\bm{p}}M \times T_{\bm{p}}M \to N_{\bm{p}}M
$$
is defined by
$$
II_{\bm p} = \langle\bm\nu_1,d^2 X\rangle\bm\nu_1 + \langle\bm\nu_2,d^2 X\rangle\bm\nu_2,
$$
where 
$$
\langle\bm\nu_i,d^2 X\rangle = l_i\,du^2+2m_i\,du\,dv+n_i\,dv^2 \quad (i=1,2).
$$

The \textit{normal curvature vector} $\kappa_n$ at $\bm{p}$ is defined by
$$
\kappa_n(\bm{p},\bm{w}) = \dfrac{II_{\bm{p}}(\bm{w},\bm{w})}{I_{\bm{p}}(\bm{w},\bm{w})},
$$
where $\bm{w} \in T_{\bm{p}}M\setminus\{\bm{0}\}$. 
The image of the unit circle in $T_{\bm{p}}M$ under $\kappa_n(\bm{p},\bm{w})$ forms an ellipse in $N_{\bm{p}}M$, which may degenerate to a point or a radial line segment. 
This image is called the \textit{curvature ellipse} at $\bm{p}$ and denoted by $\mathcal{E}_{\bm{p}}$. 
Hence, 
\begin{equation*}
\mathcal{E}_{\bm{p}}=\{\kappa_n(\bm{p},\bm{w}): \bm{w}\in T_{\bm{p}}M, I_{\bm{p}}(\bm{w},\bm{w})=1\}.
\end{equation*}

Points $\bm{p}\in M$ are classified by the position of the origin $\bm{0}_{\bm{p}}$ of $N_{\bm{p}}M$ with respect to a non-degenerate curvature ellipse $\mathcal{E}_{\bm{p}}$.  
The point $\bm{p}$ is called the \textit{elliptic}, \textit{parabolic} or \textit{hyperbolic} point if $\bm{0}_{\bm{p}}$ lies inside, on or outside $\mathcal{E}_{\bm{p}}$, respectively. 
On the other hand, $\bm{p}$ is called an \textit{inflection} point if $\mathcal{E}_{\bm{p}}$ degenerates to a radial line segment or a point.   

The curvature ellipse $\mathcal{E}_{\bm{p}}$ can be viewed as the image of the unit circle in $T_{\bm{p}}M$ under a pair of quadratic forms $(Q_1,Q_2)$ associated with $II_{\bm{p}}$.
Writing $\bm{w}=u X_(\bm{q})+v X_v(\bm{q})$, we show that a pair of binary quadratic forms is determined by $II_{\bm{p}}$ as follows:
$$
Q_i(u,v)=l_i u^2+2m_i u v+n_i v^2 \quad (i=1,2).
$$ 
To study the geometry of $M$, the pencil of binary forms determined by $(Q_1,Q_2)$ was investigated in \cite{BN1998}.
Identifying a binary form $A u^2+2B u v+C v^2$ with its coefficients $(A,B,C)\in\R^3$, we see that there is a double cone defined by the equation $B^2-AC=0$, representing the perfect squares (up to sign), in $\R^3$.
This cone determines a conic in $\R P^2$. 
If $Q_1$ and $Q_2$ are linearly independent, then their pencil determines a projective line in $\R P^2$, 
and thus $\bm{p}$ is referred to as elliptic, parabolic or hyperbolic point if this line intersects the conic in zero, one or two points, respectively, 
On the other hand, if $Q_1$ and $Q_2$ are linearly dependent, then $\bm{p}$ is referred to as inflection point if the rank of $\mathcal A(\bm{p})$ is at most one (see \cite{Little1969}), where 
\begin{equation*}
\mathcal{A}(\bm{p}) = 
\begin{pmatrix}
l_1 & m_1 & n_1 \\
l_2 & m_2 & n_2
\end{pmatrix}
\end{equation*}
We remark that if $(Q_1,Q_2)$ is equivalent to $(x^2,0)$ or $(0,0)$ under the action of $\mathcal{G}=GL(2,\R) \times GL(2,\R)$, then $\bm{p}$ is called a \textit{degenerate inflection} point (see \cite{BT2002, NBT2007}). 
We also remark that if $(Q_1,Q_2)$ is $\mathcal{G}$-equivalent to $(x^2,0)$, then $\mathcal{E}_{\bm{p}}$ degenerates to a radial line segment and $\bm{0}_{\bm{p}}$ is an endpoint of $\mathcal{E}_{\bm{p}}$, whereas if $(Q_1, Q_2)$ is $\mathcal{G}$-equivalent to $(0,0)$, then $\mathcal{E}_{\bm{p}}$ degenerates to a point $\bm{0}_{\bm{p}}$. 
In this case, $\bm{p}$ is an inflection point in the sense of \cite{Little1969}, but not in the sense of \cite{MRFR1995}. 

Let $\bm{w}(\theta)\in S(T_{\bm{p}}M)$, $\theta\in[0,2\pi)$, be a regular parameterization of the unit circle in $T_{\bm{p}}M$. 
A direction $\bm{w}(\theta_0) \in T_{\bm{p}}M$ along which 
$$
\kappa_n(\bm{p},\bm{w}(\theta_0))\quad \text{and} \quad \dfrac{\partial}{\partial \theta}\kappa_n(\bm{p},\bm{w}(\theta_0))
$$
are parallel is an \textit{asymptotic direction} of $M$ at $\bm{p}$, where $\bm{w}(\theta)$ and $-\bm{w}(\theta)$ are identified as the same direction (see \cite{GMRFR2000, MRFR1995}).  

For a unit normal vector $\bm\nu \in S(N_{\bm{p}}M)$, a symmetric bilinear form 
$$
II_{\bm{p},\bm\nu} \colon T_{\bm{p}}M \times T_{\bm{p}}M \to \R
$$
is defined by
$$
II_{\bm{p},\bm\nu}(\bm{w}_1,\bm{w}_2)=\langle II_{\bm{p}}(\bm{w}_1,\bm{w}_2), \bm\nu \rangle,\quad (\bm{w}_i \in T_{\bm{p}}M).
$$
Hence, in terms of the tangent frame $\{X_u,X_v\}$ of $T_{\bm{p}}M$, we  see that
\begin{equation}
\label{eq:II_pv}
\begin{split}
\left[II_{\bm{p},\bm\nu}\right] &= 
\begin{pmatrix}
II_{\bm{p},\bm\nu}(X_u,X_u) & II_{\bm{p},\bm\nu}(X_u,X_v) \\
II_{\bm{p},\bm\nu}(X_u,X_v) & II_{\bm{p},\bm\nu}(X_v,X_v)
\end{pmatrix}
=
\begin{pmatrix}
\langle X_{uu}, \bm\nu \rangle & \langle X_{uv}, \bm\nu \rangle \\
\langle X_{uv}, \bm\nu \rangle & \langle X_{vv}, \bm\nu \rangle
\end{pmatrix}.
\end{split}
\end{equation}
A normal direction $\bm\nu$ is called a \textit{binormal direction} of $M$ at $\bm{p}$ if the height function $h_{\bm\nu}(u,v)=\langle X(u,v), \bm\nu \rangle$ on $M$ along $\bm\nu$ has a degenerate singularity at $\bm{q}$.  
This characterization is also holds at inflection points (see \cite{MRFR1995}) and at points where the curvature ellipse degenerates to a point (which are inflection points in the sense of \cite{Little1969}), even if a unique asymptotic direction cannot be determined.

\subsection{Contact}
\label{sec:contact}

We recall the notion of contact, in terms of singularity theory, introduced by Montaldi \cite{Montaldi1986-1}. 
For a more detailed and specific treatment, we refer the reader to \cite[Chapter~4]{IRFRT2015}.

Let $(X_1,Y_1)$ and $(X_2,Y_2)$ be two pairs of submanifold-germs at the origin of $\R^n$.
The pairs are said to {\it have the same contact type} if there exists a diffeomorphism-germ $\alpha \colon (\R^n, \bm{0}) \to (\R^n, \bm{0})$ such that $\alpha(X_1)=X_2$ and $\alpha(Y_1)=Y_2$.

The following shows that contact type is measured by $\mathcal{K}$-equivalence. 
Here, two map-germs $f_1, f_2 \colon (\R^n, \bm{0}) \to (\R^p,\bm{0})$ are {\it $\mathcal{K}$-equivalent} if there exist a diffeomorphism-germ $\alpha \colon (\R^n, \bm{0}) \to (\R^n, \bm{0})$ and a matrix valued map $A \colon (\R^n, \bm{0}) \to GL(p,\R)$ such that $A f_1=f_2\circ \alpha.$

\begin{theorem}[\cite{Montaldi1986-1}]
\label{thm:contact}
Let $\iota_i \colon (X_i,\bm{x}_i) \to (\R^n, \bm{0})$ $(i=1,2)$ be immersion-germs and let $\beta_i\colon (\R^n, \bm{0}) \to (\R^p, \bm{0})$ be submersion-germs with $Y_i=\beta_i^{-1}(\bm{0})$.
Then the pairs $(X_1,Y_1)$ and $(X_2,Y_2)$ have the same contact type if and only if $\beta_1 \circ \iota_1$ and $\beta_2 \circ \iota_2$ are $\mathcal{K}$-equivalent.
\end{theorem}

By virtue of Theorem~\ref{thm:contact}, if $\iota \colon (X,\bm{x}) \to (\R^n,\bm{0})$ is an immersion-germ and $\beta \colon (\R^n, \bm{0}) \to (\R^p, \bm{0})$ is a submersion-germ with $Y=\beta^{-1}(\bm{0})$, then the contact type of $X$ and $Y$ is determined by the $\mathcal{K}$-equivalence class of the map-germ $\beta \circ \iota$. 
We call the map-germ the {\it contact map-germ} of $X$ and $Y$.

In this paper, we consider the following $\mathcal{K}$-singularities (see \cite{Arnold1986}):
$$
\setlength{\extrarowheight}{6pt} 
\begin{array}{ 
>{\centering\arraybackslash$}m{2.4cm}<{$} | >{\centering\arraybackslash$}m{2.4cm}<{$} | >{\centering\arraybackslash$}m{2.4cm}<{$} | >{\centering\arraybackslash$}m{2.4cm}<{$} | >{\centering\arraybackslash$}m{2.4cm}<{$} }
A_k^{\pm}, \quad k\geq1 & D_k^{\pm}, \quad k \geq 4 & E_6 & E_7 & E_8 \\ \hline 
x^2 \pm y^{k+1} & x^2y \pm y^{k-1} & x^3 \pm y^4 & x^3+x y^3 & x^3+y^5
\end{array}
$$


\subsection{Cylindrical directions of regular surfaces in $\R^3$}
\label{sec:cylindrical_direction}

In \cite{FHN2017}, the authors introduced the cylindrical directions of regular surfaces $M_3$ in $\R^3$ via their contact with cylinders.

Let $C^{(3)}_{\bm{v},\bm{a},r}$ denote a cylinder whose axis is parallel to $\bm{v} \in S^2$, passes through $\bm{a} \in \bm{v}^{\perp}$ but does not pass through $\bm{p} \in M_3$ and whose cross-sectional circle has radius $|r| \ne0$. 

Assume that $\bm{p}$ is a common point of $M_3$ and $C^{(3)}_{\bm{v},\bm{a},r}$. 
Then $A_{\geq1}$-contact at $\bm{p}$ if and only if $M_3$ and $C^{(3)}_{\bm{v},\bm{a},r}$ have a common normal line at $\bm{p}$. 
When this is the case, $\bm{a}$ lies in the projected line of the common normal line there under the orthogonal projection along $\bm{v}$ onto $\bm{v}^{\perp}$. 
Moreover, in this case, $\bm{v} \in T_{\bm{p}} M_3$, and hence the following theorems hold:

\begin{theorem}[{\cite[Theorem 4.1]{FHN2017}}]
\label{thm:A_2}
The contact type of $M_3$ with $C^{(3)}_{\bm{v},\bm{a},r}$ at $\bm{p}$ is $A_{\geq 2}$-contact if and only if one of the following conditions holds: 
\begin{enumerate}
\item
$\bm{p}$ is not a parabolic point, $\bm{v}$ is not an asymptotic direction at $\bm{p}$ and $\bm{a}$ is the focal point of the apparent contour of $M_3$ along $\bm{v}$ at a point corresponding to $\bm{p}$. 

\item
$\bm{p}$ is a parabolic point but not an umbilic point, $\bm{v}$ is an asymptotic direction at $\bm{p}$ and $\bm{a}$ is not a point corresponding to the focal point at $\bm{p}$. 

\item
$\bm{p}$ is a flat umbilic point. 
\end{enumerate}
\end{theorem}

\begin{theorem}[{\cite[Proposition 4.12]{FHN2017}}]
\label{thm:A_3}
If Theorem~\ref{thm:A_2}~(1) holds, then the contact type of $M_3$ with $C^{(3)}_{\bm{v},\bm{a},r}$ at $\bm{p}$ is $A_{\geq k}$-contact $(k \geq 3)$ if and only if the apparent contour has a $(k-2)$-order vertex at a point corresponding to $\bm{p}$.
\end{theorem}

The direction $\bm{v}$ of the axis of $C^{(3)}_{\bm{v},\bm{a},r}$ is a \textit{cylindrical direction} of $M_3$ at $\bm{p}$ if Theorem~\ref{thm:A_2} (1) holds and the contact between $M_3$ and  $C^{(3)}_{\bm{v},\bm{a},r}$ at $\bm{p}$ is of type $A_{\geq 3}$.
More specifically, this definition requires that the apparent contour of $M_3$ along $\bm{v}$ be a regular curve and have no inflection point at the point corresponding to $\bm{p}$, that this point be a vertex of the apparent contour, and that $\bm{a}$ is its finite focal point.


\section{Contact of regular surfaces in $\R^4$ with toroidal cylinders}
\label{sec:contact_with_toroidal_cylinders}
The set of 2-planes through a fixed point in $\R^4$ is a 4-dimensional manifold, which is the Grassmannian $Gr(2,4)$. 
Nogueira \cite{Nogueira1998} used a scalar product on $\R^4$ to obtain a local parameterization of $Gr(2,4)$ and analyzed the relationship between the contact of surfaces in $\R^4$ with 2-planes and the singularity of parallel projections along planes using a family of the projections.  
While, Deolindo-Silva and Salarinoghabi \cite{DSS2026} used the affine structure of $\R^4$ to characterize singularities of the parallel projections and to describe geometric information about the surface both when the planes, which determine the directions of the projections, contain an asymptotic direction and when these do not. 
They also studied the apparent contours of the surfaces along the planes through the second fundamental form of the surface given by the canonical projection on the normal space. 
In contrast to these studies, we use a local parameterization of the oriented Grassmannian $Gr^+(2,4)\cong SO(4)/(SO(2)\times SO(2))$ to obtain a positively oriented orthonormal frame of the orthogonal complement $\pi^\perp$ of a plane $\pi$. 
By using this frame, we define the toroidal-cylindrical contact map-germ and the orthogonal projection along $\pi$, preserving the Euclidean geometric structure relevant to the contact of surfaces in $\R^4$ with toroidal cylinders and the geometry of their apparent contours.

In \cite{FHN2017}, to study the contact of a regular surface  in $\R^3$ with cylinders, the contact map-germ defined by the difference between the third components of Monge forms of the surface and the cylinder was used.  
In contrast, following the concept of Montaldi's contact map-germ, we use the contact map-germ obtained by composing a parameterization of a regular surface in $\R^4$ with a defining function whose zero set is the toroidal cylinder. 

Assume that $X \colon (\R^2,\bm{q}) \to (\R^4,\bm{p})$ is a smooth map-germ parameterizing a regular surface $M$ in $\R^4$. 
Let $\pi$ and $\pi^\perp$ be planes so that $\R^4=\pi \oplus \pi^\perp$, and let $\{\bm{v}_1,\bm{v}_2\}$ denote a positively oriented orthonormal frame of $\pi^\perp$. 
Let $C_{\pi,\bm{a},r}$ denote a toroidal cylinder whose radius is $|r|\ne0$ and axis is parallel to the plane $\pi$ and passes through a point $\bm{a}\in \pi^\perp$ but does not pass through $\bm{p}$. 
We then define $C_{\pi,\bm{a},r}$ by 
$$
C_{\pi,\bm{a},r}=\{\bm{x}\in\R^4\colon \lVert \langle \bm{x}, \bm{v}_1\rangle \bm{v}_1 + \langle \bm{x},\bm{v}_2\rangle \bm{v}_2 -\bm{a} \rVert^2=r^2\}.
$$
Let $\bm{p} \in M \cap C_{\pi,\bm{a},r}$. 
The contact map-germ of $M$ and $C_{\pi,\bm{a},r}$ is given by 
$$
F_{\pi,\bm{a}}(u,v)=\lVert \langle X(u,v), \bm{v}_1 \rangle \bm{v}_1+\langle X(u,v),\bm{v}_2\rangle \bm{v}_2-\bm{a}\rVert^2,
$$
which we call the \textit{toroidal-cylindrical contact map-germ}.
Since toroidal cylinders are determined by their axes $\pi$, which are elements of the Grassmannian $Gr(2,4)$, reference points $\bm{a}\in \pi^\perp$ through which their axes pass and radii, the moduli space of toroidal cylinders is seven-dimensional. 
Therefore, we expect that there are $A_{\leq 7}$, $D_{\leq 7}$ and $E_{\leq 8}$-contact toroidal cylinders with $M$. 

Let $f_{\pi}$ denote the parallel projection of $M$ onto $\pi^\perp$ defined by 
$$
f_{\pi} =P_{\pi} \circ X \colon (\R^2,\bm{q}) \to \pi^\perp,\quad f_{\pi}(u,v) = \langle X(u,v),\bm{v}_1 \rangle \bm{v}_1 + \langle X(u,v), \bm{v}_2 \rangle \bm{v}_2,
$$
where $P_{\pi}$ denotes the parallel projection onto $\pi^\perp$: 
$$
P_{\pi} \colon \R^4 \to \pi^\perp, \quad P_{\pi}(\bm{x})=\langle \bm{x},\bm{v}_1 \rangle \bm{v}_1+\langle \bm{x},\bm{v}_2\rangle \bm{v}_2.
$$
Note that, by the definitions of $F_{\pi,\bm{a}}$ and $f_\pi$, we see that
$$
F_{\pi,\bm{a}}(u,v)=\lVert f_\pi(u,v) - \bm{a} \rVert^2.
$$

Locally, parallel projections can be regarded as map-germs from $(\R^2,\bm{0}) \to (\R^2,\bm{0})$, whose singularities are classified under the $\mathcal{A}$-equivalence. 
Here, two map-germs $f_1, f_2\colon(\R^2,\bm{0})\to(\R^2,\bm{0})$ are said to be \textit{$\mathcal{A}$-equivalent} if $f_2=\Psi \circ f_1 \circ \Phi^{-1}$ for some diffeomorphisms-germs $\Phi$ and $\Psi$ of the source and target, respectively. 
Rieger \cite{Rieger1987} gave a classification of corank-one map-germs (Table~\ref{tab:corank1}), while Rieger and Ruas \cite{RR1991} gave a classification of corank-two map-germs (Table~\ref{tab:corank2}). 

\begin{center}
\begin{table}[ht!]
\caption{$\mathcal{A}$-singularities of $(\R^2,\bm{0}) \to (\R^2,\bm{0})$ of corank one with $\mathcal{A}_e$-codimension less than or equal to four (\cite{Rieger1987}).}
\setlength{\extrarowheight}{6pt}
\begin{tabular}{ccc} 
Type & Normal form & $\mathcal{A}_e$-codimension \\ \hline
$1$ (Immersion) & $(x,y)$ & $0$ \\
$2$ (Fold) & $(x,y^2)$ & $0$ \\
$3$ (Cusp) & $(x,xy+y^3)$ & $0$ \\
$4_2$ (Lipse/Beaks) & $(x,y^3\pm x^2y)$ & $1$ \\
$5$ (Swallowtail) & $(x,xy+y^4)$ & $1$ \\
$4_3$ (Goose) & $(x,y^3+x^3y)$ & $2$ \\
$6$ (Butterfly) & $(x,xy+y^5\pm y^7)$ & $2$ \\
$11_5$ (Gulls) & $(x,xy^2+y^4+y^5)$ & $2$ \\
$4_4$ (Ugly goose) & $(x,y^3\pm x^4y)$ & $3$ \\
$7$ (Elder butterfly) & $(x,xy+y^5)$ & $3$ \\ 
$8$ & $(x,xy+y^6\pm y^8+\alpha y^9)$ & $3^\dagger$ \\ 
$11_7$ (Ugly gulls) & $(x,xy^2+y^4+y^7)$ & $3$ \\
$12$ & $(x,xy^2+y^5+y^6)$ & $3$ \\
$16$ & $(x,x^2y+y^4\pm y^5)$ & $3$ \\ 
$4_5$ & $(x,y^3+x^5y)$ & $4$ \\
$9$ & $(x,xy+y^6+y^9)$ & $4$ \\
$10^\ast$ & $(x,xy+y^7\pm y^9+\alpha y^{10}+\beta y^{11})$ & $4^\dagger$ \\ 
$11_9$ & $(x,xy^2+y^4+y^9)$ & $4$ \\
$13$ & $(x,xy^2+y^5 \pm y^9)$ & $4$ \\
$15$ & $(x,xy^2+y^6+y^7+\alpha y^9)$ & $4^\dagger$ \\
$17$ & $(x,x^2y+y^4)$ & $4$ \\
$18^\ast$ & $(x,x^2y+xy^3+\alpha y^5+y^6+\beta y^7)$ & $4^\dagger$ \\
$19$& $(x,x^3y+\alpha x^2 y^2+y^4+x^3 y^2)$ & $4^\dagger$ \\ \hline
\end{tabular}\\
\vspace*{3pt}
${}^\dagger$codimension of is that of $\mathcal{A}$-stratum \\
${}^\ast$excluding exceptional values of the moduli
\label{tab:corank1}
\end{table}
\end{center}

\begin{center}
\begin{table}[ht!]
\caption{$\mathcal{A}$-singularities of $(\R^2,\bm{0}) \to (\R^2,\bm{0})$ of corank two with $\mathcal{A}_e$-codimension less than or equal to four (\cite{RR1991}).}
\setlength{\extrarowheight}{6pt}
\begin{tabular}{ccc} 
Type & Normal form & $\mathcal{A}_e$-codimension \\ \hline
$I^{1,1}_{2,2}$ (Sharksfin) & $(x^2+y^3,y^2+x^3)$ & $2$ \\
$I^{2,1}_{2,2}$ (Odd-shaped sharksfin) & $(x^2+y^5,y^2+x^3)$ & $3$ \\
$I^{3,1}_{2,2}$ & $(x^2+y^7,y^2+x^3)$ & $4$ \\
$I^{2,2}_{2,2}$ & $(x^2+y^5,y^2+x^5)$ & $4$ \\
$II^1_{2,2}$ (Deltoid) & $(x^2-y^2+x^3,xy)$ & $2$ \\ 
$II^2_{2,2}$ & $(x^2-y^2+x^5,xy)$ & $4$ \\ \hline
\end{tabular}
\label{tab:corank2}
\end{table}
\end{center}

Firstly, we give a criterion for the singularities of $F_{\pi,\bm{a}}$ divided into the cases $\text{rank}(df_\pi(\bm{q}))$ is one and zero. 

\begin{proposition}
\label{prop:A1}
The contact map-germ $F_{\pi,\bm{a}}$ has an $A_{\geq 1}$-singularity at $\bm{q}$  if and only if one of the following holds: 
\begin{enumerate}
\item
$\mathrm{rank}(d f_{\pi}(\bm{q}))=1$ and there exists a unit vector $\bm\nu \in \pi^\perp \cap L^\perp$ such that 
$$
\bm{a}=f_{\pi}(\bm{q})+r\bm\nu,
$$
where $L^\perp$ denotes the orthogonal complement of $L=\mathrm{Im}(d f_{\pi}(\bm{q}))$.

\item
$\mathrm{rank}(d f_{\pi}(\bm{q}))=0$ and there exists a unit vector $\bm\nu \in \pi^\perp=N_{\bm{p}}M$ such that 
$$
\bm{a}=f_{\pi}(\bm{q})+r\bm\nu.
$$ 
\end{enumerate}
\end{proposition}
\begin{proof}
Since $\{\bm{v}_1,\bm{v}_2\}$ is a positively oriented orthonormal frame  of $\pi^\perp$, we obtain  
\begin{align}
\label{eq:1-dri}
\begin{split}
(F_{\pi,\bm{a}})_u &= 2 \left( \langle X_u, \bm{v}_1 \rangle \langle X-\bm{a}, \bm{v}_1 \rangle + \langle X_u, \bm{v}_2 \rangle \langle X-\bm{a}, \bm{v}_2 \rangle \right), \\
(F_{\pi,\bm{a}})_v &= 2 \left( \langle X_v, \bm{v}_1 \rangle \langle X-\bm{a}, \bm{v}_1 \rangle + \langle X_v, \bm{v}_2 \rangle \langle X-\bm{a}, \bm{v}_2 \rangle\right). 
\end{split}
\end{align}
Therefore, $F_{\pi,\bm{a}}$ has an $A_{\geq1}$-singularity at $\bm{q}$ if and only if 
\begin{equation}
\label{eq:dF}
\begin{pmatrix}
\langle X_u(\bm{q}),\bm{v}_1 \rangle & \langle X_u(\bm{q}),\bm{v}_2 \rangle \\
\langle X_v(\bm{q}),\bm{v}_1 \rangle & \langle X_v(\bm{q}),\bm{v}_2 \rangle
\end{pmatrix}
\begin{pmatrix}
\langle X(\bm{q})-\bm{a},\bm{v}_1 \rangle \\
\langle X(\bm{q})-\bm{a},\bm{v}_2 \rangle
\end{pmatrix}
=
\begin{pmatrix}
0 \\ 0
\end{pmatrix}.
\end{equation}
Since $|r|\ne0$, \eqref{eq:dF} implies  
\begin{equation}
\label{eq:A1}
\begin{vmatrix}
\langle X_u(\bm{q}),\bm{v}_1 \rangle & \langle X_u(\bm{q}),\bm{v}_2 \rangle \\
\langle X_v(\bm{q}),\bm{v}_1 \rangle & \langle X_v(\bm{q}),\bm{v}_2 \rangle
\end{vmatrix}
=0. 
\end{equation}

Using the orthonormal frame defined above, we identify $\pi^\perp$ with $\R^2$ and the parallel projection $f_{\pi}$ can be also expressed by 
$$
f_{\pi}(u,v)=(\langle X(u,v), \bm{v}_1 \rangle, \langle X(u,v), \bm{v}_2 \rangle ).
$$
Hence, $df_{\pi}(\bm{q})$ is given by 
$$
df_{\pi}(\bm{q}) = 
\begin{pmatrix}
\langle X_u(\bm{q}),\bm{v}_1 \rangle & \langle X_u(\bm{q}),\bm{v}_2 \rangle \\
\langle X_v(\bm{q}),\bm{v}_1 \rangle & \langle X_v(\bm{q}),\bm{v}_2 \rangle
\end{pmatrix}^T,
$$
and thus $f_{\pi}$ has a singularity at $\bm{q}$ if and only if \eqref{eq:A1} holds. 

The equation \eqref{eq:dF} is equivalent to 
\begin{equation}
\label{eq:A1-2}
\begin{pmatrix}
\langle (f_{\pi})_u(\bm{q}),\bm{v}_1 \rangle & \langle (f_{\pi})_u(\bm{q}),\bm{v}_2 \rangle \\
\langle (f_{\pi})_v(\bm{q}),\bm{v}_1 \rangle & \langle (f_{\pi})_v(\bm{q}),\bm{v}_2 \rangle
\end{pmatrix}
\begin{pmatrix}
\langle f_{\pi}(\bm{q})-\bm{a},\bm{v}_1 \rangle \\
\langle f_{\pi}(\bm{q})-\bm{a},\bm{v}_2 \rangle
\end{pmatrix}
=
\begin{pmatrix}
0 \\ 0
\end{pmatrix}.
\end{equation}
Assume that $\mathrm{rank}(df_{\pi}(\bm{q}))=1$. 
If 
$$
(\langle (f_{\pi})_u(\bm{q}), \bm{v}_1 \rangle, \langle (f_{\pi})_u(\bm{q}), \bm{v}_2 \rangle)\ne(0,0),
$$ 
then, by the definition of $f_\pi$,
\begin{align*}
&\langle X_u(\bm{q}),\bm{v}_1\rangle \langle X(\bm{q})-\bm{a},\bm{v}_1\rangle + \langle X_u(\bm{q}),\bm{v}_2\rangle
\langle X(\bm{q})-\bm{a},\bm{v}_2\rangle \\
&=\langle\langle X_u(\bm{q}),\bm{v}_1\rangle \bm{v}_1+\langle X_u(\bm{q}),\bm{v}_2\rangle \bm{v}_2, X(\bm{q})-\bm{a} \rangle =\langle (f_{\pi})_u(\bm{q}), f_{\pi}(\bm{q})-\bm{a}\rangle = 0.
\end{align*}
Since $(f_{\pi})_u(\bm{q})$ and $(f_{\pi})_v(\bm{q})$ are linearly dependent and $\bm{p}\in C_{\pi,\bm{a},r}$, equation \eqref{eq:A1-2} holds if and only if
\begin{equation}
\label{eq:a}
\bm{a} = f_{\pi}(\bm{q})+r\bm\nu
\end{equation}
for a suitable unit vector $\bm\nu \in \pi^\perp \cap L^\perp$. 
The same conclusion holds when 
$$
(\langle (f_{\pi})_v(\bm{q}), \bm{v}_1 \rangle, \langle (f_{\pi})_v(\bm{q}), \bm{v}_2 \rangle)\ne(0,0).
$$ 
We prove assertion (1).

Assume that $\mathrm{rank}(d f_{\pi}(\bm{q}))=0$. 
Since 
$$
(f_\pi)_u(\bm{q})=(f_\pi)_v(\bm{q})=0, 
$$
\eqref{eq:A1-2} holds automatically. 
Moreover, since $\bm{p}\in C_{\pi,\bm{a},r}$, there exists a suitable unit vector $\bm\nu \in \pi^\perp$ such that \eqref{eq:a} holds. 
Furthermore, it follows from the regularity of $X$ that if $\text{rank}(df_\pi(\bm{q}))=0$, then $\pi=T_{\bm{p}}M$, and thus $\pi^\perp=N_{\bm{p}}M$.  
This completes the proof. 
\end{proof}

\begin{theorem}
\label{thm:A2}
Assume that $F_{\pi,\bm{a}}$ has an $A_{\geq 1}$-singularity at $\bm{q}$. 
Then $F_{\pi,\bm{a}}$ has an $A_{\geq2}$-singularity at $\bm{q}$ if and only if one of the following holds:
\begin{enumerate}
\item
$\mathrm{rank}(df_{\pi}(\bm{q})) = 1$, $f_{\pi}$ has a fold singularity at $\bm{q}$ and 
\begin{equation}
\label{eq:A2_a}
\bm{a} = f_\pi(\bm{q})+\frac1{\kappa(\bm{q})}\bm\nu,
\end{equation}
where $\kappa$ denotes the curvature of the apparent contour. 

\item
$\mathrm{rank}(df_{\pi}(\bm{q})) = 1$, $\bm{q}$ is a singular point of the Jacobian determinant $J(u,v) = \det(df_{\pi}(u,v))$ and either
$$
\bm{a} \ne f_\pi(\bm{q}) + \dfrac1{\kappa_{\pi,\bm\nu}(\bm{q})}\bm\nu \quad (\kappa_{\pi,\bm\nu}(\bm{q})\ne0)
$$
or $\mathcal{E}_{\bm{p}}$ is included in the line which is perpendicular to $\bm\nu$ and passes through $\bm{0}_{\bm{p}}$. 
Here, the definition of $\kappa_{\pi,\bm\nu}$ is given in the proof.

\item
$\mathrm{rank}(df_{\pi}(\bm{q})) = 0$ and $\bm\nu \in N_{\bm{p}}M$ is a binormal direction of $M$ at a parabolic or hyperbolic point of $M$ at $\bm{p}$.

\item
$\mathrm{rank}(df_{\pi}(\bm{q})) = 0$, $(Q_1,Q_2)$ is equivalent to $(x^2,0)$ under the action $\mathcal{G}$, and $\mathrm{rank}(H(h_{\bm\nu})(\bm{q}))=1$.
\end{enumerate}
\end{theorem}

\begin{proof}
Differentiating \eqref{eq:1-dri} once more, we obtain
\begin{align*}
(F_{\pi,\bm{a}})_{uu} &= 2 \sum_{i=1}^2 \left( \langle X_{uu}, \bm{v}_i \rangle \langle X - \bm{a}, \bm{v}_i \rangle + \langle X_u, \bm{v}_i \rangle^2 \right), \\
(F_{\pi,\bm{a}})_{uv} &= 2 \sum_{i=1}^2 \left( \langle X_{uv}, \bm{v}_i \rangle \langle X - \bm{a}, \bm{v}_i \rangle + \langle X_u, \bm{v}_i \rangle \langle X_v, \bm{v}_i \rangle \right), \\
(F_{\pi,\bm{a}})_{vv} &= 2 \sum_{i=1}^2 \left( \langle X_{vv}, \bm{v}_i \rangle \langle X - \bm{a}, \bm{v}_i \rangle + \langle X_v, \bm{v}_i \rangle^2 \right).
\end{align*}
Since $\bm{a} \in \pi^\perp$ and $\{\bm{v}_1,\bm{v}_2\}$ is an orthonormal frame of $\pi^\perp$, using $f_{\pi} = \langle X, \bm{v}_1 \rangle \bm{v}_1 + \langle X, \bm{v}_2 \rangle \bm{v}_2$, the Hessian matrix of $F_{\pi,\bm{a}}$ at $\bm{q}$ can be written as
\begin{align*}
H(F_{\pi,\bm{a}}) = 2
\begin{pmatrix}
\langle (f_{\pi})_{uu}, f_{\pi} - \bm{a} \rangle + \lVert (f_{\pi})_u \rVert^2 &
\langle (f_{\pi})_{uv}, f_{\pi} - \bm{a} \rangle + \langle (f_{\pi})_u, (f_{\pi})_v \rangle \\
\langle (f_{\pi})_{uv}, f_{\pi} - \bm{a} \rangle + \langle (f_{\pi})_u, (f_{\pi})_v \rangle &
\langle (f_{\pi})_{vv}, f_{\pi} - \bm{a} \rangle + \lVert (f_{\pi})_v \rVert^2
\end{pmatrix}.
\end{align*}
Substituting $f_{\pi}(\bm{q}) - \bm{a} = -r\bm{\nu}$ into the above matrix yields
\begin{equation}
\label{eq:hess}
H(F_{\pi,\bm{a}}) = 2
\begin{pmatrix}
\lVert (f_{\pi})_u \rVert^2 - r\langle (f_{\pi})_{uu}, \bm{\nu} \rangle &
\langle (f_{\pi})_u, (f_{\pi})_v \rangle - r\langle (f_{\pi})_{uv}, \bm{\nu} \rangle \\
\langle (f_{\pi})_u, (f_{\pi})_v \rangle - r\langle (f_{\pi})_{uv}, \bm{\nu} \rangle &
\lVert (f_{\pi})_v \rVert^2 - r\langle (f_{\pi})_{vv}, \bm{\nu} \rangle
\end{pmatrix}.
\end{equation}
Note that $F_{\pi,\bm{a}}$ has a degenerate singularity (i.e., non-Morse singularity) at $\bm{q}$ if and only if $\det(H(F_{\pi,\bm{a}})(\bm{q})) = 0$.

First, assume that $\mathrm{rank}(df_{\pi}(\bm{q})) = 1$. 
By an appropriate change of local coordinates $(u,v)$ around $\bm{q}$, we may assume without loss of generality that $(f_{\pi})_u(\bm{q}) = \bm{0}$ and $(f_{\pi})_v(\bm{q}) \neq \bm{0}$. 
Then \eqref{eq:hess} reduces to
\begin{equation}
\label{eq:hess-2}
H(F_{\pi,\bm{a}})(\bm{q}) = 2
\begin{pmatrix}
-r\langle (f_{\pi})_{uu}(\bm{q}), \bm{\nu} \rangle &
-r\langle (f_{\pi})_{uv}(\bm{q}), \bm{\nu} \rangle \\
-r\langle (f_{\pi})_{uv}(\bm{q}), \bm{\nu} \rangle &
\lVert (f_{\pi})_v(\bm{q}) \rVert^2 - r\langle (f_{\pi})_{vv}(\bm{q}), \bm{\nu} \rangle
\end{pmatrix}.
\end{equation}
Since $|r| \ne 0$, the condition $\det(H(F_{\pi,\bm{a}})(\bm{q})) = 0$ is equivalent to
\begin{equation}
\label{eq:A2_rank1}
r(P R-Q^2)-P\lVert (f_\pi)_v(\bm{q})\rVert^2=0,
\end{equation}
where 
$$
P=\langle (f_{\pi})_{uu}(\bm{q}), \bm{\nu} \rangle, \quad Q= \langle (f_{\pi})_{uv}(\bm{q}), \bm{\nu} \rangle, \quad  R=\langle (f_{\pi})_{vv}(\bm{q}), \bm{\nu} \rangle.
$$

Let $J(u,v) = \det(df_{\pi}(u,v))$ be the Jacobian determinant of $f_{\pi}$ with respect to the orthonormal frame $\{\bm{v}_1,\bm{v}_2\}$. 
Since $(f_{\pi})_u(\bm{q}) = \bm{0}$, we have
$$
J_u(\bm{q}) = \det((f_{\pi})_{uu}(\bm{q}), (f_{\pi})_v(\bm{q})), \quad
J_v(\bm{q}) = \det((f_{\pi})_{uv}(\bm{q}), (f_{\pi})_v(\bm{q})).
$$
Since $\bm{\nu} \in \pi^\perp$ is a unit vector orthogonal to $L = \text{Im}(df_{\pi}(\bm{q})) = \text{span}\{(f_{\pi})_v(\bm{q})\}$, it follows that $J_u(\bm{q}) = 0$ if and only if $P = 0$, and similarly $J_v(\bm{q}) = 0$ if and only if $Q = 0$.

Since $(f_{\pi})_u(\bm{q}) = \bm{0}$ and $(f_{\pi})_v(\bm{q}) \ne \bm{0}$, the apparent contour of $M$ along $\pi$ is regular at $f_\pi(\bm{q})$ (i.e., the singularity of $f_\pi$ is a fold at $\bm{q}$ if and only if $P\ne0$. 
In this case, the singular set $\Sigma(f_\pi)=J^{-1}(0)$ of $f_\pi$ is parameterized near $\bm{q}$ by $(u(v),v)$ with $u'(\bm{q})=-J_v(\bm{q})/J_u(\bm{q})$ by the implicit function theorem. 
Hence, the apparent contour $\gamma(v)=f_\pi(u(v),v)$ near $f_\pi(\bm{q})$ satisfying 
\begin{align*}
\gamma'(v) &= u'(v)(f_{\pi})_u(u(v),v) + (f_{\pi})_v(u(v),v), \\
\gamma''(v) &= u''(v)(f_{\pi})_u(u(v),v) + (u'(v))^2(f_{\pi})_{uu}(u(v),v) \\
&\qquad+ 2u'(v)(f_{\pi})_{uv}(u(v),v) + (f_{\pi})_{vv}(u(v),v)
\end{align*}
Evaluating at $\bm{q}$ where $(f_{\pi})_u(\bm{q}) = \bm{0}$, we have 
$$
\gamma'(\bm{q}) = (f_{\pi})_v(\bm{q}), \quad \gamma''(\bm{q}) = (u'(\bm{q}))^2(f_{\pi})_{uu}(\bm{q}) + 2u'(\bm{q})(f_{\pi})_{uv}(\bm{q}) + (f_{\pi})_{vv}(\bm{q}).
$$
Using 
$$
u'(\bm{q}) = -\frac{J_v(\bm{q})}{J_u(\bm{q})} = -\frac{\langle (f_{\pi})_{uv}(\bm{q}), \bm{\nu} \rangle}{\langle (f_{\pi})_{uu}(\bm{q}), \bm{\nu} \rangle}, 
$$
the curvature $\kappa(\bm{q})$ of $\gamma$ with respect to the unit normal vector $\bm\nu$ is given by
$$
\kappa(\bm{q})=\frac{\langle \gamma''(\bm{q}), \bm{\nu} \rangle}{\lVert \gamma'(\bm{q}) \rVert^2}
= \frac{PR-Q^2}{P \lVert (f_{\pi})_v(\bm{q}) \rVert^2}.
$$
If $\gamma$ does not have an inflection point at $f_\pi(\bm{q})$, that is, $PR-Q^2\ne0$, then substituting this into \eqref{eq:A2_rank1}, we obtain $r\,\kappa(\bm{q}) - 1 = 0$. 
Since $|r| \neq 0$, this holds if and only if $\kappa(\bm{q}) \neq 0$ and $r = 1/\kappa(\bm{q})$. 
In this case, $\bm{a} = f_{\pi}(\bm{q}) + 1/{\kappa(\bm{q})}\bm{\nu}$ is the focal point of the apparent contour at $f_{\pi}(\bm{q})$. 
In contrast, if $\gamma$ has an inflection point $f_\pi(\bm{q})$, then the left hand side of \eqref{eq:A2_rank1} reduces to $-P \lVert (f_\pi)_v(\bm{q}) \rVert^2$. 
Since $P\ne0$, \eqref{eq:A2_rank1} does not hold.    

The Jacobian determinant $J$ is regular at $\bm{q}$ but the apparent contour $f_\pi(\Sigma(f_\pi))$ is singular at $f_\pi(\bm{q})$ if and only if $P=0$ and $Q\ne0$. 
In this case, the left hand side of \eqref{eq:A2_rank1} reduces $-r Q^2$. 
Since $|r|\ne0$, \eqref{eq:A2_rank1} does not hold.

Next, we assume that $J$ has a singularity at $\bm{q}$, that is, $dJ(\bm{q})=0$.  
If follows form \eqref{eq:II_pv} and \eqref{eq:hess} that 
\begin{equation}
\label{eq:hess_3}
H(F_{\pi,\bm{a}}(\bm{q})) = 2\left(\left[\tilde{I}_{\bm{p}}\right]-r\left[ II_{\bm{p},\bm\nu}\right]\right),
\end{equation}
where $\tilde{I}_{\bm{p}}$ denotes the first fundamental form according to $f_{\pi}$ which is defined by 
$$
\tilde{I}_{\bm{p}}=(df_\pi(\bm{q}))^T df_\pi(\bm{q})\colon T_{\bm{p}}M \times T_{\bm{p}}M \to \R. 
$$ 
Via $dX_{\bm{q}}$, we identify $\R^2$ with $T_{\bm{p}}M$ by 
$$
(u,v) \leftrightarrow u X_u(\bm{q})+v X_v(\bm{q}).
$$ 
Since $\tilde{I}_{\bm{p}}$ is a pullback of the Euclidean metric by $f_\pi$, 
$$
\tilde{I}_{\bm{p}}(\bm\eta, \bm{w}) = \langle df_\pi(\bm\eta),df_\pi(\bm{w})\rangle,
$$
where $\bm\eta \in \text{ker}(df_\pi(\bm{q}))\subset T_{\bm{p}}M$, and 
$$
\text{ker}(\tilde{I}_{\bm{p}}) = \text{ker}(df_\pi(\bm{q})), \quad \text{rank}(\tilde{I}_{\bm{p}})=1. 
$$ 
Moreover, since $dJ(\bm{q})=0$, by the preceding relations between $dJ(\bm{q})$ and the coefficients $P$, $Q$, we have
$$
II_{\bm{p},\bm\nu}(\bm\eta,\bm{w})=0.
$$
Hence, 
$$
\text{ker}(\tilde{I}_{\bm{p}}) \subset \text{ker}(II_{\bm{p},\bm\nu}).
$$
It follows that there exists $c\in \R$ such that 
$$
II_{\bm{p},\bm\nu}=c\tilde{I}_{\bm{p}},  
$$
and thus, for every $\bm{w}\in T_{\bm{p}}M$ with $df_\pi(\bm{w})\ne0$, 
$$
c=\dfrac{{II}_{\bm{p},\bm\nu}(\bm{w},\bm{w})}{\tilde{I}_{\bm{p}}(\bm{w},\bm{w})}=\kappa_{\pi,\bm\nu}(\bm{p}).
$$ 
Substituting this into \eqref{eq:hess_3}, we obtain
$$
H(F_{\pi,\bm{a}}(\bm{q})) = 2(1-r \kappa_{\pi,\bm\nu}(\bm{q}))\left[\tilde{I}_{\bm{p}}\right].
$$
Therefore, since $\text{rank}(\tilde{I}_{\bm{p}})=1$, we see that $\text{rank}(H(F_{\pi,\bm{a}}(\bm{q})))=1$ if and only if $1/r\ne \kappa_{\pi,\bm\nu}(\bm{q})$, whereas $\text{rank}(H(F_{\pi,\bm{a}}(\bm{q}))=0$ if and only if $1/r = \kappa_{\pi,\bm\nu}(\bm{q})$. 
In particular, if $II_{\bm{p},\bm\nu}(\bm{w},\bm{w})=0$, then $\text{rank}(H(F_{\pi,\bm{a}}(\bm{q})))=1$. 
Note that for every $\bm{w}\in T_{\bm{p}}M$, $II_{\bm{p},\bm\nu}(\bm{w},\bm{w})=0$ is equivalent to $\langle II_{\bm{p}}(\bm{w},\bm{w}),\bm\nu \rangle =0$. 
In this case, $\mathcal{E}_{\bm{p}} \subset N_{\bm{p}}M \cap \bm\nu^\perp$ which is equivalent to $\mathcal{E}_{\bm{p}}$ is included in the line which is perpendicular to $\bm\nu$ and passes through $\bm{0}_{\bm{p}}$.  


Finally, assume that $\text{rank}(df_{\pi}(\bm{q})) = 0$. 
Then $(f_{\pi})_u(\bm{q}) = (f_{\pi})_v(\bm{q}) = \bm{0}$ and $\pi^\perp = N_{\bm{p}}M$. 
Since $\bm{\nu} \in \pi^\perp$, we have $\langle (f_{\pi})_{uu}(\bm{q}), \bm{\nu} \rangle = \langle X_{uu}(\bm{q}), \bm{\nu} \rangle$, and similarly for the other second-order derivatives. 
Thus, \eqref{eq:hess} becomes
$$
H(F_{\pi,\bm{a}})(\bm{q}) = -2r
\begin{pmatrix}
\langle X_{uu}(\bm{q}), \bm{\nu} \rangle & \langle X_{uv}(\bm{q}), \bm{\nu} \rangle \\
\langle X_{uv}(\bm{q}), \bm{\nu} \rangle & \langle X_{vv}(\bm{q}), \bm{\nu} \rangle
\end{pmatrix},
$$
which is $-2r$ times the Hessian matrix $H(h_{\bm\nu})$ of the height function $h_{\bm{\nu}}(u,v) = \langle X(u,v), \bm{\nu} \rangle$ on $M$ along $\bm\nu$ at $\bm{q}$. 
Since $|r| \ne 0$, $\det(H(F_{\pi,\bm{a}})(\bm{q})) = 0$ holds if and only if $\det(H(h_{\bm{\nu}})(\bm{q})) = 0$, meaning that $\bm{\nu} \in N_{\bm{p}}M$ is a binormal direction of $M$ at $\bm{p}$. 
For a generic surface, $(Q_1, Q_2)$ is equivalent to neither $(x^2,0)$ nor $(0,0)$ under the action $\mathcal{G}$. 
In this context, $h_{\bm\nu}$ has an $A_{\geq2}$-singularity at $\bm{q}$ when $\bm{q}$ is a parabolic or hyperbolic point of $M$, whereas $h_{\bm\nu}$ has a $D_{\geq4}$- or more degenerate singularity at $\bm{q}$ (i.e., $\text{rank}(H(h_{\bm\nu})(\bm{q}))=0$) if and only if $Q_1$ and $Q_2$ are dependent (see \cite{BN1998}).  
When $(Q_1,Q_2)$ is equivalent to $(0,0)$, since the second fundamental form vanishes, $\text{rank}(H(h_{\bm\nu})(\bm{q}))=0$ for every $\bm\nu \in N_{\bm{p}}M$. 
On the other hand, when $(Q_1,Q_2)$ is equivalent to $(x^2,0)$, $\text{rank}(H(h_{\bm\nu})(\bm{q}))=1$ for a suitable $\bm\nu$.

\end{proof}

\begin{remark}
\label{rem:A2}
It follows from the proof of Theorem~\ref{thm:A2} that if $f_{\pi}$ has a fold singularity at $\bm{q}$ whose apparent contour has an inflection point at $f_{\pi}(\bm{q})$, or if $J$ is regular at $\bm{q}$ but the apparent contour $f_\pi(\Sigma(f_\pi))$ is singular at $f_\pi(\bm{q})$, then $F_{\pi,\bm{a}}$ cannot have an $A_{\geq 2}$-singularity at $\bm{q}$ for any $|r| \neq 0$; in these cases, $F_{\pi,\bm{a}}$ always has an $A_1$-singularity at $\bm{q}$. 
In contrast to case the contact of regular surfaces in $\R^3$ with cylinders, these are distinctive properties (see \cite[Theorem~4.1]{FHN2017} for details).
\end{remark}

The following theorem shows that when an apparent contour of $M$ along $\pi$ is both regular and non-inflectional, the singularity of $F_{\pi,\bm{a}}$ is characterized by the order of the vertex of the apparent contour. 
The key idea of the proof is similar to those in \cite[Proposition~4.12]{FHN2017} and \cite[Proposition~2.6]{KRSUY2005}. 

\begin{theorem}
\label{thm:Ak}
Assume that Theorem~\ref{thm:A2} (1) holds.
The function $F_{\pi,\bm{a}}$ has an $A_{k \geq3}$-singularity at $\bm{q}$ if and only if the apparent contour of $M$ along $\pi$ has a $(k-2)$-order vertex at $f_\pi(\bm{q})$.
\end{theorem}

\begin{proof}
Since $f_\pi$ has a fold singularity at $\bm{q}$, its null direction is transverse to the singular set $\Sigma(f_\pi)$ of $f_\pi$. 
Hence, by choosing local coordinates $(s,t)$ with $(0,0)=\bm{q}$ as in the proof of \cite[Proposition~2.6]{KRSUY2005}, we may assume that $\Sigma(f_\pi)=\{t=0\}$ and $\text{ker}(df_\pi)$ along $\Sigma(f_\pi)$ is spanned by $\{\partial_t\}$. 
We see that $(f_\pi)_t(s,0)=0$, and thus we obtain $(F_{\pi,\bm{a}})_t(s,0)=0$. 
Moreover, since $f_\pi$ is a fold and \eqref{eq:A2_a} holds, 
$$
(F_{\pi,\bm{a}})_{tt}(0,0)=-2\langle (f_\pi)_{tt}(0,0),\bm{a}-f_\pi(0,0) \rangle =-\dfrac2{\kappa(0,0)}\langle (f_\pi)_{tt}(0,0),\bm\nu \rangle \ne0.
$$ 
Therefore, by the parametric Morse lemma, we obtain 
\begin{equation}
\label{eq:Fs}
F_{\pi,\bm{a}}(s,t)\sim_{\mathcal{R}} \pm t^2 +F_{\pi,\bm{a}}(s,0)
\end{equation}
under a suitable change of coordinate. 
Here, two function-germ $g_1, g_2\colon (\R^2,\bm{0})\to(\R,0)$ are said to be \textit{$\mathcal{R}$-equivalent} if $g_2=g_1\circ \Psi$ for some differomorphism-germ $\Psi \colon (\R^2,\bm{0})\to(\R^2,\bm{0})$ of the source, and we denote $g_1\sim_{\mathcal{R}} g_2$. 
Moreover, since $\Sigma(f_\pi)=\{t=0\}$, we obtain 
$$
F_{\pi,\bm{a}}(s,0)=\lVert f_\pi(s,0) -\bm{a} \rVert^2.
$$
It follows from this equation and \eqref{eq:Fs} that 
$$
F_{\pi,\bm{a}}(s,t)\sim_\mathcal{R} \pm t^2+F_{\pi,\bm{a}}(s,0)= \pm t^2+\lVert \gamma(s)-\bm{a} \rVert^2
$$
where $\gamma$ denotes a local parameterization of an apparent contour $f_\pi(\Sigma(f_\pi))$ along $\pi$ and $\bm{a}$ is regarded as a point $\pi^\perp \subset \R^4$, and thus $F_{\pi,\bm{a}}$ has an $A_{k\geq3}$-singularity at $\bm{q}$ if and only if $\gamma$ has a $(k+1)$-point contact with its osculating circle at $f_\pi(\bm{q})$. 
In this case, $\gamma$ has a $(k-2)$-order vertex at $f_\pi(\bm{q})$ (see, for example, \cite{BG1992}). 
\end{proof}

Theorem~\ref{thm:Ak} yields that the axis, which is an affine plane $\bm{a}+\pi$, of an $A_{\geq3}$-contact toroidal cylinder $C_{\pi,\bm{a},r}$ of a regular surface in $\R^4$ is an analogue of the axis of an $A_{\geq 3}$-contact cylinder of a regular surface in $\R^3$. 
Therefore, the notion of cylindrical directions of regular surfaces in $\R^3$ can be extended to regular surfaces in $\R^4$ as follows: 
Assume that $f_\pi$ is fold at $\bm{q}$. 
Then the direction $\pi$ of the axis of the $A_{\geq3}$-contact toroidal cylinder is called the \textit{toroidal-cylindrical direction} of a regular surface in $\R^4$ at $\bm{p}$ derived from its contact with the toroidal cylinder $C_{\pi,\bm{a},r}$. 
  
\begin{remark}
As shown in Remark~\ref{rem:A2}, if $f_\pi$ has a fold singularity at $\bm{q}$ and its apparent contour has an inflection point at $f_\pi(\bm{q})$, then the contact between $M$ and $C_{\pi,\bm{a},r}$ is of type $A_1$. 
Therefore, under the fold condition, $A_{\geq2}$-contact automatically implies that the apparent contour is non-inflectional and that $\bm{a}$ is its finite focal point. 
Consequently, neither condition needs to be imposed separately in the definition of a toroidal-cylindrical direction. 
This contrasts with the case of cylindrical directions of regular surfaces in $\R^3$.
\end{remark}

\end{document}